\documentclass{amsart}
\usepackage[utf8]{inputenc}
\usepackage{amsmath}
\usepackage{amsfonts}
\usepackage{amssymb}
\usepackage{amsthm}
\usepackage{graphicx}
\usepackage[colorinlistoftodos]{todonotes}
\usepackage[colorlinks=true, allcolors=blue]{hyperref}
\usepackage{enumerate}
\usepackage{bm}
\usepackage{enumitem}

\newtheorem{thm}{Theorem}[section]
\newtheorem{lem}[thm]{Lemma}
\newtheorem{cor}[thm]{Corollary}

\newtheorem{prop}[thm]{Proposition}
\newtheorem{rmk}[thm]{Remark}

\makeatletter
\newcommand*{\bdiv}{%
  \nonscript\mskip-\medmuskip\mkern5mu%
  \mathbin{\operator@font div}\penalty900\mkern5mu%
  \nonscript\mskip-\medmuskip
}
\makeatother

\theoremstyle{definition}

\title[Existence and Morse Index of two free boundary embedded geodesics]{Existence and Morse Index of two free boundary embedded geodesics on Riemannian 2-disks with convex boundary}
\author{Dongyeong Ko}
\address{Department of Mathematics, Rutgers University - New Brunswick, Piscataway, NJ 08854}
\email{dk954@math.rutgers.edu}
\begin{document}

\maketitle
\begin{abstract}
We prove that a free boundary curve shortening flow on closed surfaces with a strictly convex boundary remains noncollapsed for a finite time in the sense of the reflected chord-arc profile introduced by Langford-Zhu. This shows that such flow converges to free boundary embedded geodesic in infinite time, or shrinks to a round half-point on the boundary. As a consequence, we prove the existence of two free boundary embedded geodesics on a Riemannian $2$-disk with a strictly convex boundary. Moreover, we prove that there exists a simple closed geodesic with Morse Index $1$ and $2$. This settles the free boundary analog of Grayson's theorem.
\end{abstract}
\section{Introduction}

The search of closed geodesics on surfaces by min-max methods was initiated by Birkhoff \cite{B}. The following classical theorem on the existence of three embedded geodesics on spheres was initially proven by Lusternik-Schnirelmann \cite{LS}, and the embeddedness part is proven by Grayson \cite{Gr} by curve shortening flow.

\begin{thm} [Lusternik-Schnirelmann \cite{LS}, Grayson \cite{Gr}]
A closed Riemannian $2$-sphere $(S^{2},g)$ contains at least three simple closed geodesics.
\end{thm}

In the free boundary setting, Lusternik-Schnirelmann \cite{LS} proved that a bounded domain in $\mathbb{R}^{n}$ with a convex and smooth boundary admits at least $n$ distinct geodesic chords which meet the boundary orthogonally. Bos \cite{Bo} extended the existence of orthogonal geodesic chords to Riemannian manifolds via Birkhoff's discrete curve shortening process. He also proved the necessity of the convexity assumption by giving an example of a non-convex domain in $\mathbb{R}^{2}$ (See Figure 1 in \cite{Bo}). The construction of free boundary (immersed) geodesics in various settings have been developed by Gluck-Ziller \cite{GZ} and Zhou \cite{Z} by the discrete curve shortening process. Li \cite{L} introduced the chord shortening flow and provided a simpler approach to prove Lusternik-Schnirelmann's existence result. Moreover, the recent work of Donato-Montezuma \cite{DM} constructed a free boundary embedded geodesic or a geodesic loop in the nonnegative sectional curvature scenario via Almgren-Pitts min-max theory (See also \cite{AMS}). However, to the author's knowledge, the existence of free boundary embedded geodesics on Riemannian $2$-disks with convex boundary is not known. The following is our first main result on this: 

\begin{thm}
There are at least $2$ free boundary embedded geodesics on a Riemannian $2$-disk $(D^{2},\partial D^{2},g)$ with a strictly convex boundary.
\end{thm}
For $a>b>0$, since the ellipse $E := \{ (x,y) | x^{2}/a^{2} + y^{2}/b^{2}=1 \}$ achieves only two free boundary geodesics, the existence theorem of two free boundary embedded geodesics (Theorem 1.1) is optimal.

A Riemannian metric $g$ is called \emph{bumpy} if every closed geodesic is nondegenerate i.e. there is no free boundary geodesic that admits a non-trivial Jacobi field. The arguments in Abraham \cite{Ab} and Ambrozio-Carlotto-Sharp \cite{ACS} shows that the set of bumpy metrics is generic in $C^{r}$-sense for $r \ge 5$. We have
\begin{cor}
There are at least two free boundary embedded geodesics with distinct length on the Riemannian $2$-disk endowed with a bumpy metric and strictly convex boundary.
\end{cor}

We note here that we do not need to impose the condition ensures the nonexistence of interior simple closed geodesics because the mass cancellation such as neck-pinching over the limiting process does not happen in the curve setting. This is different from the min-max construction of free boundary minimal disks in Haslhofer-Ketover \cite{HK}. 

We compare our result with the problem on its higher dimensional analogue. There is a conjecture on the existence of three free boundary minimal disks on Riemannian $3$-balls with strictly convex boundary and nonnegative Ricci curvature. Struwe \cite{St} constructed one free boundary immersed minimal disk via mapping approach. Grüter and Jost \cite{GJ} constructed a free boundary minimal surface by developing Simon-Smith min-max theory in the free boundary setting (See also Jost \cite{J}). We refer to \cite{DR}, \cite{F}, \cite{LP}, \cite{Li2}, \cite{LiZ}, \cite{LSZ} for seminal works on the construction of free boundary minimal disks. Recently, Haslhofer-Ketover \cite{HK} proved this conjecture for generic metric and obtained at least $2$ free boundary minimal disks by combining the ideas from min-max theory, free boundary mean curvature flow and degree theory. 

To obtain the existence of a free boundary embedded geodesic, we apply the free boundary curve shortening flow for the tightening procedure. The study of free boundary curve shortening flow (FBCSF), the Neumann boundary problem of curve shortening flow, dates back to the work of Huisken \cite{H1} on free boundary mean curvature flow on the graph setting in general dimensions. After this work, free boundary mean curvature flow is extensively studied in \cite{AW}, \cite{St1} and \cite{St2}. One of the most recent development of free boundary curve shortening flow is the full determination of the long time behavior of the free boundary curve shortening flow on convex domain in $\mathbb{R}^{2}$ of Langford-Zhu \cite{LZ}. We generalize this result to the surfaces with strictly convex boundary. 

Huisken \cite{H2} found a simpler proof of the long-time behavior of curve shortening flows for closed and embedded curves by developing the idea of distance comparison principle. The essence of their idea is to prove that the ratio between extrinsic and intrinsic distance has monotonicity so that does not collapse along the flow (See also \cite{AB}). Edelen \cite{E} extended Huisken's work to the surface setting. He proved that the ratio along the flow decreases at worst exponentially while the ratio may not have the monotonicity. Langford-Zhu \cite{LZ} developed the distance comparison principle in the free boundary setting by introducing the notion of reflected chord-arc profile. Our generalization of the idea of reflected chord-arc profile to the surface setting proves that the chord-arc profile achieves the lower bound weighted by the exponential terms by time over the lower bound in the Euclidean setting, where our approach relies on Edelen's work. The noncollapsing property gives the following:

\begin{thm} Let $(N,\partial N, g)$ be a closed Riemannian surface with strictly convex boundary and $\{ \Gamma_{t} \}_{t \in [0,T)}$ be a maximal free boundary curve shortening flow starting from a properly embedded closed interval $\Gamma_{0}$ in $N$. Then either:
\begin{enumerate}
    \item $T = \infty$, in which case $\Gamma_{t}$ converges smoothly as $t \rightarrow \infty$ to an embedded geodesic in $N$ which meets $\partial N$ orthogonally; or
    \item $T< \infty$, in which case $\Gamma_{t}$ converges uniformly to some single round half-point $z \in \partial N$ smoothly in the sense of the blow up limit of the curve converges to the unit semi-circle.
\end{enumerate}
\end{thm}

Marques-Neves (\cite{MN1}, \cite{MN2}) proved the upper and lower bound of Morse Index of min-max minimal hypersurfaces on closed manifolds with dimension $3 \le n+1 \le 7$ by proving deformation theorems. In the smooth setting of curves, we are forced to construct the interpolation deformation between two families of simple closed curves with controlled length, where the deformation arises from the curve shortening flow. The author \cite{K1} proved the Morse Index bound of the simple closed geodesics on Riemannian $2$-spheres (see also \cite{DMMS}). Moreover, the author \cite{K2} proved the Morse Index bound of min-max capillary embedded geodesics on certain Riemannian $2$-disks by constructing the interpolation based on the curve shortening flow with a strictly fixed boundary. 

\begin{thm}
Suppose $(D^{2},\partial D^{2}, g)$ is a Riemannian $2$-disk with convex boundary endowed with a bumpy metric. Then for each $k=1,2$, there exists a free boundary embedded geodesic $\gamma_{k}$ with
\begin{equation*}
    index(\gamma_{k}) = k
\end{equation*}
and these two geodesics satisfy $|\gamma_{1}| <|\gamma_{2}|$.
\end{thm}
\begin{cor}
For a $2$-Riemannian disk $(D^{2},\partial D^{2},g)$ with a convex boundary, for $k=1,2$, there exists a free boundary embedded geodesic $\gamma_{k}$ with
\begin{equation*}
    index(\gamma_{k}) \le k \le index(\gamma_{k})+ nullity (\gamma_{k}).
\end{equation*}
\end{cor}

By the work of Smale \cite{Sm} on the diffeomorphisms of $2$-sphere, the space of embedded intervals on $D^{2}$ after identifying the point curves retract onto $\mathbb{R}P^{2}$. Lusternik-Schnirelmann arguments, for instance in \cite{CLOT}, suggests that the number of critical points of a smooth real-valued function has the lower bound of one plus the maximal cup-length. In our setting, the cohomology ring of $\mathbb{R}P^{2}$ has a maximal length $2$, there are at least $3$ critical points. The point curves can be regarded as a stable critical point and we expect there are at least $2$ free boundary embedded geodesics.

We prove the existence of two free boundary embedded geodesics by considering one and two parameter min-max construction. The free boundary curve shortening flow provides the way to tight the minimizing sequences and the classical Lusternik-Schnirelmann arguments gives the existence. 

Our interpolation arguments rely on the quantitative estimate of $F$-distance bound in \cite{K1} along the squeezing homotopy arising from the free boundary curve shortening flow. The main modification we needed to make was the construction of free boundary mean convex foliation when the geodesic is strictly stable. We obtain a free boundary mean convex foliation by the first eigenfunction of the stability operator with Robin boundary condition (See also Proposition 2.4 in \cite{HK}).

The organization of the paper is as follows. In Section 2, we collect the basic notions on free boundary curve shortening flow, the reflected chord-arc profile and the second variation of a free boundary embedded geodesic. In Section 3, we discuss the variation of the chord-arc profile on surfaces. In Section 4, we prove the noncollapsedness of the free boundary curve shortening flow at finite times. In Section 5, we consider the min-max construction of free boundary embedded geodesics. In Section 6, we prove the Morse Index bound of free boundary embedded geodesics.

\section*{Acknowledgments}
The author wishes to thank his advisor Daniel Ketover for his constant support and valuable discussions. The author also thanks to Otis Chodosh for helpful explanations relating to this work. The author also thanks to Jonathan Zhu for valuable discussion relating to this work. The author was partially supported by NSF grant DMS-1906385. 
\section{Preliminaries}
\subsection{Free boundary curve shortening flow}
We consider an oriented Riemannian surface with boundary $(N^{2},\partial N^{2},g)$ with a convex $C^{2}$ boundary and a properly immersed family of curves with boundary $ \gamma: M^{1} \times I \rightarrow N^{2}$. We denote $J$ by the counterclockwise rotation by $\pi/2$ and take the convention that the orientation of the unit normal vector field on $\nu$ is given to satisfy $\gamma'/|\gamma'|=J\nu$. We call a properly immersed family of curves with boundary $\{ \Gamma_{t} \}$ satisfies a \emph{free boundary curve shortening flow} if $\gamma: M^{1} \times I  \rightarrow N^{2}$ with $\Gamma_{t} = \gamma(M,t)$ satisfies
\begin{align*}
    \begin{cases}
    \partial_{t} \gamma = \kappa \nu \text{ in } \mathring{M} \times I \\
    \langle \nu, \mathcal{N}^{\partial N} \rangle = 0 \text{ on } \partial M \times I,
    \end{cases}
\end{align*}
where $\kappa(\cdot,t)$ is a geodesic curvature of $\Gamma_{t}$ with respect to the unit normal vector field $\nu$, and $\mathcal{N}^{\partial N}$ is the inward unit normal vector field on $\partial N$. We define $\kappa^{\partial N}$ to be a geodesic curvature of boundary $\partial N$. We consider the setting that $\gamma(\cdot,t)$ are embeddings and have two boundary points on $\partial N$. Since $\partial N$ is convex, there is no loss of generality to assume that $\Gamma_{t}$ does not have touching points in the interior of the curve on the boundary by the maximum principle.

We recall some notions related to \emph{completed chord-arc profile} in \cite{LZ}. For a given (unit length) parametrization of curve $\gamma: M \rightarrow N$ and $x,y \in \gamma$, denote $d(x,y)$ and $l(x,y)$ to be a distance in $N$ and the arclength between $x$ and $y$, respectively. We define the \emph{reflected distance} $\tilde{d}(x,y)$ and \emph{reflected arclength} $\tilde{l}(x,y)$ between two points in $x, y \in M$ by
\begin{equation*}
    \tilde{d}(x,y) = \min_{z \in \partial N} (d(x,z)+d(y,z))
\end{equation*}
and 
\begin{equation*}
    \tilde{l}(x,y) = \min_{s \in \partial \gamma} (l(x,s)+l(y,s)).
\end{equation*}
The \emph{reflected chord-arc profile} $\tilde{\psi}_{\gamma}$ of $\gamma$ is denoted by
\begin{equation*}
    \tilde{\psi}_{\gamma}(\delta) = \inf \{ \tilde{d}(x,y) : x,y \in \gamma, \tilde{l}(x,y) = \delta \},
\end{equation*}
and the extended chord-arc profile $\bm{\psi}_{\gamma}$ is defined by
\begin{equation*}
    \bm{\psi}_{\gamma}(\delta) = \min \{ \psi_{\gamma}(\delta),\tilde{\psi}_{\gamma}(\delta) \}.
\end{equation*}
We consider a connected, properly immersed curve-with-boundary $\gamma$ in $N$ whose endpoints are on $\partial N$. Also in case of we need the doubling of $N$, then we denote the doubling of $N$ by $\tilde{N}$. We denote the formal double $\bm{M} = (M \sqcup M)/ \partial M$ and write $\bm{x} = (x,sign(\bm{x})) \in \bm{M}$ where $sign(\bm{x})$ distinguishes to which copy of $M$ it belongs. We also denote continuous curve $\bm{\gamma} : \bm{M} \rightarrow N$ by $\bm{\gamma}(\bm{x}) = \gamma(x)$. Then we define the \emph{completed arclength} by
\begin{equation*}
    \bm{l}(\bm{x},\bm{y}) = \begin{cases} l(\gamma(x),\gamma(y)), \text{ if } sign(\bm{x}) = sign(\bm{y}) \\ \tilde{l}(\gamma(x),\gamma(y)), \text{ if } sign(\bm{x}) \neq sign(\bm{y}).
    \end{cases}
\end{equation*}
Also we define the \emph{completed distance} function $\bm{d}(\bm{x},\bm{y})$ on $\bm{M} \times \bm{M}$ by
\begin{equation*}
    \bm{d}(\bm{x},\bm{y}) = \begin{cases} d(\gamma(x),\gamma(y)), \text{ if } sign(\bm{x}) = sign(\bm{y}) \\ \tilde{d}(\gamma(x),\gamma(y)), \text{ if } sign(\bm{x}) \neq sign(\bm{y}).
    \end{cases}
\end{equation*}
Then we now define the \emph{completed chord-arc profile} $\bm{\psi}$ of $\Gamma$ by
\begin{equation*}
    \bm{\psi}(\delta) = \inf \{ \bm{d}(\bm{x},\bm{y}) : \bm{x},\bm{y} \in \bm{M}, \bm{l}(\bm{x},\bm{y}) = \delta \}.
\end{equation*}
Also we denote $\bm{\psi}(\delta,t)$ by the completed chord-arc profile of $\Gamma_{t} = \gamma(M,t)$. We call the completed chord-arc profile $\bm{\psi}$ as \emph{a classical profile} if $sign(\bm{x}) = sign(\bm{y})$ and \emph{a reflected profile} otherwise. We define $L = |\gamma|$ and $\bm{L} = 2|\gamma|$, and we denote the length functions depending on time $t$ by $L(t) = |\Gamma_{t}|$ and $\bm{L}(t) = 2|\Gamma_{t}|$. 

We also follow the setting of an auxiliary function introduced in \cite{LZ}. First, we will obtain the control of the chord-arc profile by a $C^{2}$-function $\varphi \in C^{2}([0,1])$ satisfying the following:
\begin{enumerate}[label=(\roman*)]
\item $\varphi(1-\eta) = \varphi(\eta)$ for all $\eta \in [0,1]$.
\item $0<\varphi'<1$.
\item $\varphi$ is strictly concave.
\end{enumerate}
Note that $\varphi(\bm{l}/\bm{L})$ is smooth away from the diagonal $\bm{D}$ in $\bm{M} \times \bm{M}$ and well-defined since $0<\bm{l}(\bm{x},\bm{y})/\bm{L} \le 1/2$ for $\bm{x}, \bm{y} \in \bm{M}$.

In case of considering a time dependent auxiliary function, we consider the following function; We define a $C^{2}$-function $\varphi: [0,1] \times [0,T) \rightarrow \mathbb{R}$ satisfying the following conditions for every time $t \in [0, T)$.
\begin{enumerate}[label=(\roman*)]
    \item $\varphi(1-\eta,t) = \varphi(\eta,t)$ for all $\eta \in [0,1]$.
    \item $|\partial_{\eta}\varphi(\eta,t)|<1$ for $\eta \in [0,1]$ and $t \in [0,T)$.
    \item $\varphi(\cdot,t)$ is strictly concave.
\end{enumerate}

We consider the auxiliary functions $Z$ and $\tilde{Z}$ on $M \times M$ given by
\begin{align*}
    Z(x,y) = d(\gamma(x),\gamma(y)) - \bm{L} \varphi \Big( \frac{l(\gamma(x),\gamma(y))}{\bm{L}} \Big), \\
    \tilde{Z}(x,y) = \tilde{d}(\gamma(x),\gamma(y)) - \bm{L} \varphi \Big( \frac{\tilde{l}(\gamma(x),\gamma(y))}{\bm{L}} \Big).
\end{align*}
If we consider the auxiliary function on $M \times M \times \partial N$ given by
\begin{equation*}
    \overline{Z}(x,y,z) = d(\gamma(x),z) + d(\gamma(y),z) - \bm{L} \varphi \Big( \frac{\tilde{l}(\gamma(x),\gamma(y))}{\bm{L}} \Big)
\end{equation*}
then $\tilde{Z}(x,y) = \min_{z \in \partial N} \overline{Z}(x,y,z)$. We denote our completed auxiliary function $\bm{Z}$ on $\bm{M} \times \bm{M}$ by
\begin{equation}
    \bm{Z}(\bm{x},\bm{y}) = \bm{d}(\bm{x},\bm{y}) - \bm{L} \varphi \Big( \frac{\bm{l}(\bm{x},\bm{y})}{\bm{L}} \Big) = \begin{cases}
        Z(\bm{x},\bm{y}) \text{ if } sign(\bm{x}) = sign(\bm{y}) \\ \tilde{Z}(\bm{x},\bm{y}) \text{ if } sign(\bm{x}) \neq sign(\bm{y}).
    \end{cases}
\end{equation}
We also denote $d(\cdot,\cdot,t)$, $\tilde{d}(\cdot,\cdot,t)$, $\bm{d}(\cdot,\cdot,t)$ and $l(\cdot,\cdot,t)$, $\tilde{l}(\cdot,\cdot,t)$, $\bm{l}(\cdot,\cdot,t)$ by the distances, lengths of $\Gamma_{t}$.
We consider the auxiliary functions at time $t$ as
\begin{align*}
    Z(x,y,t) = d(\gamma(x),\gamma(y),t) - \bm{L}(t) \varphi \Big( \frac{l(\gamma(x),\gamma(y),t)}{\bm{L}} \Big), \\
    \tilde{Z}(x,y,t) = \tilde{d}(\gamma(x),\gamma(y),t) - \bm{L}(t) \varphi \Big( \frac{\tilde{l}(\gamma(x),\gamma(y),t)}{\bm{L}} \Big).
\end{align*}
Then the auxiliary function $\overline{Z}$ at time $t$ is defined by
\begin{equation*}
    \overline{Z}(x,y,z,t) = d(\gamma(x),z,t) + d(\gamma(y),z,t) - \bm{L} \varphi \Big( \frac{\tilde{l}(\gamma(x),\gamma(y),t)}{\bm{L}} \Big).
\end{equation*}
Then the completed auxiliary function $\bm{Z}$ at time $t$ is
\begin{equation}
    \bm{Z}(\bm{x},\bm{y},t) = \bm{d}(\bm{x},\bm{y},t) - \bm{L} \varphi \Big( \frac{\bm{l}(\bm{x},\bm{y},t)}{\bm{L}} \Big) = \begin{cases}
        Z(x,y,t) \text{ if } sign(\bm{x}) = sign(\bm{y}) \\ \tilde{Z}(x,y,t) \text{ if } sign(\bm{x}) \neq sign(\bm{y}).
    \end{cases}
\end{equation}

\subsection{Second Variation, Morse Index and the first eigenfunction}

We consider a free boundary embedded geodesic $\gamma$ on a Riemannian $2$-disk $(D^{2},\partial D^{2},g)$. Let us denote $f$ by a smooth section of the normal bundle of $\gamma$. Then the second variation of the length of $\gamma$ is defined as
\begin{align}
    \nonumber Q(f,f) &:=  \int_{\gamma} (|\nabla_{\gamma}f|^{2} - K f^{2}) ds - (\kappa(p_{1})f^{2}(p_{1}) + \kappa(p_{2})f^{2}(p_{2})) \\ 
    &=  - \int_{\gamma} (f  Lf) ds + f(p_{1})(-\nabla_{\gamma}f(p_{1})-\kappa(p_{1})f(p_{1})) + f(p_{2})(\nabla_{\gamma}f(p_{2})-\kappa(p_{2})f(p_{2})).
\end{align}
where $K$ is a Gaussian curvature on $D^{2}$, $\kappa$ is a geodesic curvature of $\partial D^{2}$ and $L = \Delta_{\gamma} + K$ is the Jacobi operator of $\gamma$. The boundary condition of the Jacobi operator is
\begin{equation*}
    (-1)^{i}\nabla_{\gamma}f(p_{i})-\kappa(p_{i})f(p_{i})=0
\end{equation*} for $i=1,2$. Let us consider the increasing sequence of eigenvalues $\{ \lambda_{i} \}$ and associated eigenfunctions $\{ \phi_{i} \}$ of the following equation with Robin boundary condition which correspond to the eigenvalues of the stability operator $Q$:
\begin{equation}
    \begin{cases} L \phi_{i} + \lambda_{i} \phi_{i} = 0 \text{ on }\gamma\\
    \frac{\partial \phi_{i}}{\partial \eta} -\kappa \phi_{i} = 0 \text{ on } \partial \gamma,
    \end{cases}
\end{equation}
where $\eta$ is an outward unit vector on $\gamma$. We define Morse index of $\gamma$ to be a maximal dimension of a subspace of a space of $C^{\infty}(\gamma)$ where $Q$ is negative definite. i.e. the number of negative eigenvalues of the stability operator $Q$. Note that the first eigenfunction $\phi_{1}$ of (4) can be chosen to be strictly positive by the standard elliptic theory.

\section{The spatial variation of the chord-arc profile on surfaces}
In this section, we calculate the variation of the chord-arc profile on surfaces, in particular, we find inequalities which holds if the auxiliary function $\bm{Z}$ achieves a zero minimizer at some pair of points $(\bm{x}_{0},\bm{y}_{0}) \in (\bm{M} \times \bm{M}) \setminus \bm{D}$. This is a generalization of Section 4 of \cite{LZ} to surfaces. We denote $\alpha_{(\bm{x},\bm{y})}$ by the (possibly broken) geodesic realizing $\bm{d}(\bm{x},\bm{y})$ parametrized (with a constant speed) by $\alpha_{(\bm{x},\bm{y})}(0) = \bm{\gamma}(\bm{x})$ and $\alpha_{(\bm{x},\bm{y})}(1) = \bm{\gamma}(\bm{y})$ between $\bm{\gamma}(\bm{x})$ and $\bm{\gamma}(\bm{y})$. We omit the subscription $(\bm{x},\bm{y})$ for simplicity and denote this curve by $\alpha$. Let us consider the doubled curve $\bm{\alpha}$ of $\alpha$ in $\tilde{N}$ if $sign(\bm{x}) \neq sign(\bm{y})$, namely the curve connecting $\bm{\gamma}(\bm{x})$ and $z$, and $z$ and $\bm{\gamma}(\bm{y})$ in $\tilde{N}$, and $\bm{\alpha} = \alpha$ otherwise. Also denote $\partial_{s}$ by the arclength parameter on $\bm{\gamma}$. Also denote $[\bm{\gamma}(\bm{x}):\bm{\gamma}(\bm{y})]$ as the shorter portion of $\bm{\gamma} \setminus \{ \bm{\gamma}(\bm{x}), \bm{\gamma}(\bm{y}) \}$.

Before deducing the inequalities, we observe the nonexistence of the intersection point between $\bm{\alpha}$ and $[\bm{\gamma}(\bm{x}):\bm{\gamma}(\bm{y})]$ other than boundary points, which follows from the arguments in the proof of Theorem 4.1 in \cite{E}. Note that two curves $\bm{\alpha}$ and $[\bm{\gamma}(\bm{x}):\bm{\gamma}(\bm{y})]$ represent $ \bm{d}(\bm{x},\bm{y})$ and $ \bm{l}(\bm{x},\bm{y})$, respectively. We can directly apply the argument in the reflected profile cases.
\begin{prop}
There is no interior intersection point between $[\bm{\gamma}(\bm{x}):\bm{\gamma}(\bm{y})]$ and $\bm{\alpha}$.
\end{prop}
We denote the region bounded by $\bm{\alpha}$ and $[\bm{\gamma}(\bm{x}):\bm{\gamma}(\bm{y})]$ by $A_{(\bm{x},\bm{y})}$. By Proposition 3.1, $A_{(\bm{x},\bm{y})}$ is a topological disk in $\tilde{N}$. We first consider a classical profile case. We follow and modify the arguments in the proof of Theorem 4.1 of \cite{E} and Proposition 4.2 of \cite{LZ}. 

\begin{prop}
    Suppose $0 = \min _{(x,y) \in M \times M \setminus D} Z(x,y) = Z(x_{0},y_{0})$. At $(x_{0},y_{0})$, we have
    \begin{equation}
        0 \le  -4 \frac{\varphi''}{\bm{L}}- \frac{\kappa(\gamma(x_{0}))}{d} \langle \alpha'(0),\nu(\gamma(x_{0})) \rangle +\frac{\kappa(\gamma(y_{0}))}{d} \langle \alpha'(1), \nu(\gamma(y_{0}))\rangle - ( 1- \varphi'^{2} ) \int_{\alpha} K.
    \end{equation}
    \begin{proof}
    Consider a variation $\alpha_{x,y}$ of $\alpha$ which are family of curves moved by $\partial_{x}$ and $\partial_{y}$ at endpoints $\gamma(x_{0})$ and $\gamma(y_{0})$ with corresponding vector fields $V$ and $W$, respectively, satisfies the following:
    \begin{align*} 
    V(0) = -\partial_{s}(\gamma(x_{0})) \text{, } V(1) = 0 \\
    W(0) = 0 \text{, } W(1) = \partial_{s}(\gamma(y_{0})).
    \end{align*}
    Denote $\rho(\gamma(x),\gamma(y)) = |\alpha_{x,y}|$. Then by definition, $\rho(\gamma(x),\gamma(y)) \ge d(\gamma(x),\gamma(y))$. We define $\hat{Z}(x,y)$ by
    \begin{equation*}
        \hat{Z}(x,y) = \rho(\gamma(x),\gamma(y)) - \bm{L} \varphi \Big( \frac{l(\gamma(x),\gamma(y))}{\bm{L}} \Big).
    \end{equation*}
    Note that $\hat{Z} \ge Z$ for every $(x,y) \in (M \times M) \setminus D$ and $\hat{Z}(x_{0},y_{0}) = Z(x_{0},y_{0})$. Moreover, $\hat{Z}$ achieves a local minimizer and
    \begin{equation}
        \partial_{x} \hat{Z}= \partial_{y} \hat{Z} = 0
    \end{equation}
    at $(x_{0},y_{0})$. (6) yields
    \begin{equation*}
        - \partial_{x} \rho = \partial_{y} \rho = \varphi' \Big( \frac{l(\gamma(x_{0}),\gamma(y_{0}))}{\bm{L}} \Big)
    \end{equation*}
    and this implies
    \begin{equation}
\Big\langle -V(0), -\frac{\alpha'(0)}{d}\Big\rangle = \Big\langle W(1), \frac{\alpha'(1)}{d}\Big\rangle = \varphi' \Big( \frac{l(x_{0},y_{0})}{\bm{L}} \Big).  
    \end{equation}    
    Then let us consider the second variation of $\hat{Z}$ by the second variation formula of the length functional:
    \begin{align*}
         \partial_{x}^{2} \hat{Z} &= \partial_{x}^{2} \rho - \frac{1}{\bm{L}} \varphi'' \\
        &= \frac{1}{d} \Big \{ \int_{0}^{1} ( \langle (V^{\perp})', (V^{\perp})' \rangle - Rm(V^{\perp},\alpha',V^{\perp},\alpha') ) + \langle \nabla_{V}V, \alpha'\rangle|_{0}^{1}   \Big\} - \frac{1}{\bm{L}} \varphi''\\  &= \frac{1}{d} (I(V^{\perp},V^{\perp}) -\kappa(\gamma(x_{0}))  \langle \alpha'(0),\nu(\gamma(x_{0})) \rangle) - \frac{1}{\bm{L}} \varphi''.
    \end{align*}
    By the similar way, we have
    \begin{align*}
         \partial_{x} \partial_{y} \hat{Z} &= \frac{1}{d} I(V^{\perp},W^{\perp}) + \frac{1}{\bm{L}} \varphi''
         \\ \partial_{y}^{2} \hat{Z} &= \frac{1}{d} (I(W^{\perp},W^{\perp}) +\kappa(\gamma(y_{0}))  \langle \alpha'(0),\nu(\gamma(y_{0})) \rangle) - \frac{1}{\bm{L}} \varphi''.
    \end{align*}
By (7), we first have 
\begin{equation}
    |(V^{\perp} \pm W^{\perp})(0)| = |(V^{\perp} \pm W^{\perp})(1)| = \sqrt{1 - \varphi'^{2}}.
\end{equation} By Proposition 3.1, there is no interior intersection point between $\alpha$ and $[\gamma(x):\gamma(y)]$, and note that $A_{(x,y)}$ is a topological disk. Since $\gamma$ is embedded and separates $N$ into two topological disks, we see that 
\begin{equation}
    \langle \nu (\gamma(x_{0})), \alpha'(0) \rangle = \langle \nu (\gamma(y_{0})), - \alpha'(1) \rangle
\end{equation}
and so two terms in (9) have the same sign. From (8) and (9), we can set $(V^{\perp}-W^{\perp})(t)$ to be a parallel transport of $(V^{\perp}-W^{\perp})(0)$ along $\alpha$ for $t \in [0,1]$. Now Since $\hat{Z}$ achieves local minimum at $(x_{0},y_{0})$, we have
    \begin{align}
        \nonumber 0 &\le (\partial_{x} - \partial_{y})^{2} \hat{Z} \\ &= \frac{1}{d} (I(V^{\perp}-W^{\perp},V^{\perp}-W^{\perp}) -\kappa(\gamma(x_{0}))  \langle \alpha'(0),\nu(\gamma(x_{0})) \rangle +\kappa(\gamma(y_{0}))  \langle \alpha'(1),\nu(\gamma(y_{0})) \rangle) -4 \frac{\varphi''}{\bm{L}}.
    \end{align}
 Then we have
    \begin{equation}
        \frac{1}{d(x_{0},y_{0})} I(V^{\perp}-W^{\perp},V^{\perp}-W^{\perp}) = - (1- \varphi'^{2}) \int_{\alpha} K 
    \end{equation}
    since $(V^{\perp}-W^{\perp})(t)$ is a parallel transport over $\alpha$. By applying (11) into (10), we obtain (5).
    \end{proof}
\end{prop}
Now we see the case of the reflected profile. 
\begin{prop} Suppose $\bm{Z} \ge 0$ with $0 = \tilde{Z}(x_{0},y_{0})= \overline{Z}(x_{0},y_{0},z_{0})$ for some $((x_{0},y_{0}),z_{0}) \in (( \mathring{M} \times \mathring{M}) \setminus D) \times N$. Let $t_{0} =\alpha^{-1}(z_{0})$. Then at $((x_{0},y_{0}),z_{0})$,
\begin{align}
    0 \le &- (1- \varphi'^{2}) \int_{\alpha} K - \frac{\kappa(\gamma(x_{0}))}{\tilde{d}}  \langle \alpha'(0),\nu(\gamma(x_{0})) \rangle + \frac{\kappa(\gamma(y_{0}))}{\tilde{d}}  \langle \alpha'(1),\nu(\gamma(y_{0})) \rangle \\ \nonumber &- \frac{2\kappa^{\partial N}(z_{0})} {\Big\langle \frac{\alpha'_{+}(t_{0})-\alpha'_{-}(t_{0})}{\tilde{d}},\mathcal{N}^{\partial N}(z_{0}) \Big\rangle} -\frac{4}{\bm{L}} \varphi''.
\end{align}
\begin{proof}
Note that, since $\partial N$ is convex, both pieces $\alpha([0,t_{0}])$ and $\alpha([t_{0},1])$ are transversal to $\partial N$. We take the orientation of $\partial_{z}$ near $z_{0} \in \partial N$ which satisfies
\begin{equation}
    \Big\langle \partial_{z}(z_{0}), \frac{\alpha'_{-}(t_{0})}{\tilde{d}}\Big\rangle >0
\end{equation}
 on $\partial N$. Now we consider a variation $\alpha_{x,y,z}$ of $\alpha$ which are family of curves moved by $\partial_{x}$, $\partial_{y}$ and $\partial_{z}$ near $\gamma(x_{0})$, $\gamma(y_{0})$ and $z_{0}$ with the corresponding vector fields $V$, $W$ and $X$ respectively, which satisfy the following:
\begin{align*}
&V(0) = \partial_{s}(\gamma(x_{0})) \text{, } V(t) = 0 \text{ for } t \in [t_{0},1], \\  &W(t) = 0 \text{ for } t \in [0,t_{0}]\text{, }  W(1) = \partial_{s}(\gamma(y_{0})),\\ 
&X(0) = X(1) =0 \text{, } X(t_{0}) = \partial_{z}(z_{0}). 
\end{align*}
We define $\rho(x,y,z) = |\alpha_{x,y,z}|$ and define $\hat{Z}(x,y,z)$ for $(x,y,z) \in (( \mathring{M} \times \mathring{M}) \setminus D) \times N$ by
 \begin{equation*}
        \hat{Z}(x,y,z) = \rho(x,y,z) - \bm{L} \varphi \Big( \frac{\tilde{l}(\gamma(x),\gamma(y))}{\bm{L}} \Big).
    \end{equation*}
Then $\hat{Z} \ge Z$ for every $(x,y,z) \in (( \mathring{M} \times \mathring{M}) \setminus D) \times N$ and $\hat{Z}(x_{0},y_{0},z_{0}) = Z(x_{0},y_{0},z_{0})$. Moreover, $\hat{Z}(x_{0},y_{0},z_{0})$ is a local minimizer of $\hat{Z}$ and $\partial_{x}\hat{Z} = \partial_{y}\hat{Z} = \partial_{z}\hat{Z} =0$. Let us first consider the first variations. This gives
\begin{align}
    \partial_{x} \rho = \partial_{y} \rho = \varphi' \\
    \partial_{z} \rho = 0
\end{align}
and we have
\begin{align}
\Big\langle V(0), -\frac{\alpha'(0)}{\tilde{d}}\Big\rangle = \Big\langle W(1), \frac{\alpha'(1)}{\tilde{d}}\Big\rangle &= \varphi' \Big( \frac{l(x_{0},y_{0})}{\bm{L}} \Big) \\
\Big\langle X(t_{0}), \frac{\alpha'_{-}(t_{0})-\alpha'_{+}(t_{0})}{\tilde{d}}\Big\rangle &= 0.
    \end{align}
 Then by (13) and (17), there exists $\theta_{0} \in (0, \pi/2)$ such that
    \begin{equation}
        \Big\langle X(t_{0}), \frac{\alpha'_{-}(t_{0})}{\tilde{d}}\Big\rangle = \Big\langle X(t_{0}), \frac{\alpha'_{+}(t_{0})}{\tilde{d}}\Big\rangle = \cos \theta_{0}.
    \end{equation}
We calculate the second variation by $x$:
 \begin{align}
         \nonumber \partial_{x}^{2} \hat{Z} &= \partial_{x}^{2} \rho - \frac{1}{\bm{L}} \varphi'' \\
        \nonumber &= \frac{1}{\tilde{d}} \Big \{ \int_{0}^{t_{0}} ( \langle (V^{\perp})', (V^{\perp})' \rangle - Rm(V^{\perp},\alpha',V^{\perp},\alpha') ) + \langle \nabla_{V}V, \alpha'\rangle|_{0}^{t_{0}}   \Big\} - \frac{1}{\bm{L}} \varphi''\\  &= \frac{1}{\tilde{d}} (I(V^{\perp},V^{\perp}) -\kappa(\gamma(x_{0}))  \langle \alpha'(0),\nu(\gamma(x_{0})) \rangle) - \frac{1}{\bm{L}} \varphi''.
    \end{align}
    Similarly we have
    \begin{align}
        \partial_{y}^{2} \hat{Z} &= \frac{1}{\tilde{d}} (I(W^{\perp},W^{\perp}) +\kappa(\gamma(y_{0}))  \langle \alpha'(1),\nu(\gamma(y_{0})) \rangle) - \frac{1}{\bm{L}} \varphi'' \\
        \partial_{x} \partial_{y} \hat{Z} &= \frac{1}{\tilde{d}} I(V^{\perp},W^{\perp}) - \frac{1}{\bm{L}} \varphi'' \\
         \partial_{z}^{2} \hat{Z} &= \frac{1}{\tilde{d}} (I(X^{\perp},X^{\perp})+ \kappa^{\partial N}(z) \langle \alpha'_{-}(t_{0})-\alpha'_{+}(t_{0}),\mathcal{N}^{\partial N}(z_{0}) \rangle).
    \end{align}
    Also we obtain
    \begin{align}
        \nonumber \partial_{x} \partial_{z} \hat{Z} &= \partial_{x} \partial_{z} \rho \\
        \nonumber &= \frac{1}{\tilde{d}} \Big \{ \int_{0}^{t_{0}} ( \langle (V^{\perp})', (X^{\perp})' \rangle - Rm(V^{\perp},\alpha',X^{\perp},\alpha') ) + \langle \nabla_{V}X, \alpha'\rangle|_{0}^{t_{0}}   \Big\}\\  &= \frac{1}{\tilde{d}} I(V^{\perp},X^{\perp}), \\
        \partial_{y} \partial_{z} \hat{Z} &= \frac{1}{\tilde{d}} I(W^{\perp},X^{\perp}).
    \end{align}
    For the brevity of notation, we put $c = \sqrt{1- \varphi'^{2}}/ \sin \theta_{0}$. From (16) and (18), we obtain
    \begin{equation}
        |(V^{\perp} \pm W^{\perp}+cX^{\perp})(0)| = |(V^{\perp} \pm W^{\perp}+cX^{\perp})(t_{0})| = |(V^{\perp} \pm W^{\perp}+cX^{\perp})(1)| = \sqrt{1-\varphi'^{2}}.
    \end{equation}
    By applying Proposition 3.1, we can argue the same reasoning as in the proof of Proposition 3.2. Namely, we obtain
    \begin{equation}
        \langle \nu (\gamma(x_{0})), \alpha'(0) \rangle = \langle \nu (\gamma(y_{0})), - \alpha'(1) \rangle
    \end{equation}
    and we know two terms in (26) have the same sign again. Hence by (25) and (26), we can take $V, W, X$ on $[0,1]$ such that $(V^{\perp}+W^{\perp}+cX^{\perp})(t)$ to be a parallel transport of $(V^{\perp}+W^{\perp}+cX^{\perp})(0)$ and $(V^{\perp}+W^{\perp}+cX^{\perp})(t_{0})$ on $[0,t_{0})$ and $[t_{0}, 1]$, respectively. Then we have
    \begin{equation}
        \frac{1}{\tilde{d}} I(V^{\perp}+W^{\perp}+cX^{\perp},V^{\perp}+W^{\perp}+cX^{\perp}) = - (1- \varphi'^{2}) \int_{\alpha} K 
    \end{equation}

    Since $(x_{0},y_{0},z_{0})$ is a minimizer of $\hat{Z}$, we have
    \begin{align}
        \nonumber 0 &\le (\partial_{x} + \partial_{y} + c \partial_{z})^{2} \hat{Z}|_{(x_{0},y_{0},z_{0})} \\ &= \frac{1}{\tilde{d}} ( I(V^{\perp}+W^{\perp}+cX^{\perp},V^{\perp}+W^{\perp}+cX^{\perp})- \kappa(\gamma(x_{0}))  \langle \alpha'(0),\nu(\gamma(x_{0})) \rangle \\ \nonumber &+ \kappa(\gamma(y_{0}))  \langle \alpha'(1),\nu(\gamma(y_{0})) \rangle + c^{2} \kappa^{\partial N}(z) \langle \alpha'_{-}(t_{0})-\alpha'_{+}(t_{0}),\mathcal{N}^{\partial N}(z_{0}) \rangle)-\frac{4}{\bm{L}} \varphi'' \\ &= - (1- \varphi'^{2}) \int_{\alpha} K - \frac{\kappa(\gamma(x_{0}))}{\tilde{d}}  \langle \alpha'(0),\nu(\gamma(x_{0})) \rangle + \frac{\kappa(\gamma(y_{0}))}{\tilde{d}}  \langle \alpha'(1),\nu(\gamma(y_{0})) \rangle \\ \nonumber &+ \frac{c^{2} \kappa^{\partial N}(z)}{\tilde{d}} \langle \alpha'_{-}(t_{0})-\alpha'_{+}(t_{0}),\mathcal{N}^{\partial N}(z_{0}) \rangle-\frac{4}{\bm{L}} \varphi'',
    \end{align}
    where $c = \sqrt{1- \varphi^{2}}/ \sin \theta_{0}$. (28) follows by summing up (19)-(24) and (29) comes from (17). Notice that
    \begin{equation}
        \frac{1}{2}\Big\langle \frac{\alpha'_{+}(t_{0})-\alpha'_{-}(t_{0})}{\tilde{d}},\mathcal{N}^{\partial N}(z_{0}) \Big\rangle = \sin \theta_{0}.
    \end{equation}
    By applying (30) to (29), we obtain (12).
\end{proof}
\end{prop}
Now we consider the completed profile. We can still apply the arguments in \cite{LZ} to prove that the first derivatives of $\bm{Z}$ vanishes even it achieves the minimum at boundary points. We directly apply Lemma 4.4 and Lemma 4.5 in \cite{LZ} to our estimates in Proposition 3.1 and Proposition 3.2.
\begin{prop} Suppose $0 = \min_{\bm{M} \times \bm{M}} \bm{Z} = \bm{Z}(\bm{x}_{0}, \bm{y}_{0})$ for some $(\bm{x}_{0}, \bm{y}_{0}) \in (\bm{M} \times \bm{M}) \setminus \bm{D}$, then there exists $(\bm{x},\bm{y}) \in (\bm{M} \times \bm{M}) \setminus \bm{D}$ such that $\bm{Z}(\bm{x},\bm{y})=0$ and either of the following holds:
\begin{enumerate}
    \item $\text{sign}(\bm{x}) = \text{sign}(\bm{y})$ and
    \begin{equation*}
         0 \le  -4 \frac{\varphi''}{\bm{L}}- \frac{\kappa(\gamma(x))}{d} \langle \alpha'(0),\nu(\gamma(x)) \rangle +\frac{\kappa(\gamma(y))}{d} \langle \alpha'(1), \nu(\gamma(y))\rangle - ( 1- \varphi'^{2} ) \int_{\alpha} K
    \end{equation*}
    or
    \item $\text{sign}(\bm{x}) \neq \text{sign}(\bm{y})$, $\bm{x}, \bm{y} \in \mathring{\bm{M}}$, $\bm{Z}(\bm{x},\bm{y}) = \tilde{Z}(x,y)= \overline{Z}(x,y,z)$, $t_{0} = \alpha^{-1}(z)$.
    \begin{align*}
        0 \le &- (1- \varphi'^{2}) \int_{\alpha} K - \frac{\kappa(\gamma(x))}{\tilde{d}}  \langle \alpha'(0),\nu(\gamma(x)) \rangle + \frac{\kappa(\gamma(y))}{\tilde{d}}  \langle \alpha'(1),\nu(\gamma(y)) \rangle \\ &- \frac{2\kappa^{\partial N}(z)} {\Big\langle \frac{\alpha'_{+}(t_{0})-\alpha'_{-}(t_{0})}{\tilde{d}},\mathcal{N}^{\partial N}(z) \Big\rangle} -\frac{4}{\bm{L}} \varphi''.
    \end{align*}
    
\end{enumerate}
\end{prop}
\section{Noncollapsing and the long time behavior of the flow}

Based on the estimates we obtained in Proposition 3.4, we obtain the non-collapsing properties of complete chord-arc profile under the free boundary curve shortening flow on surfaces, which are generalizations of Theorem 5.3 and Theorem 5.4 in \cite{LZ}. We define $[\Gamma_{t}(x):\Gamma_{t}(y)]$ as the portion of $\Gamma_{t}$ connecting $\Gamma_{t}(x)$ and $\Gamma_{t}(y)$ for $x,y \in M$. Also denote $[\bm{\Gamma}_{t}(\bm{x}):\bm{\Gamma}_{t}(\bm{y})]$ as the shorter portion of $\bm{\Gamma}_{t} \setminus \{ \bm{\Gamma}_{t}(\bm{x}), \bm{\Gamma}_{t}(\bm{y}) \}$ as in Section 3. The same arguments with the proof of Proposition 5.1 in \cite{LZ} follows the evolution of the chord-arc profile on surfaces.

\begin{prop} 
Assume that $\bm{Z}(\cdot,\cdot,0) \ge 0$ and $\bm{Z}(\cdot,\cdot,0) >0$ on off-diagonal points. We denote $t_{0} = \sup\{ t \in [0,T):Z(\cdot,\cdot,t) \ge 0 \}< T$. Then there exist $\bm{x}, \bm{y} \in (\bm{M} \times \bm{M}) \setminus \bm{D}$ such that $\bm{Z}(\bm{x},\bm{y},t_{0}) = 0$ and either of the following holds:
\begin{enumerate}
    \item $\text{sign}(\bm{x}) = \text{sign}(\bm{y})$ and
    \begin{equation}
         0 \ge  4 \frac{\varphi''}{\bm{L}} + ( 1- \varphi'^{2} ) \int_{\alpha} K + 2 \Big( \varphi - \varphi' \frac{\bm{l}}{\bm{L}} \Big)  \int_{\Gamma_{t}} \kappa^{2} ds + \varphi' \int_{[\Gamma_{t}(x):\Gamma_{t}(y)]} \kappa^{2} ds - \bm{L} \partial_{t} \varphi
    \end{equation}
    or
    \item $\text{sign}(\bm{x}) \neq \text{sign}(\bm{y})$, $\bm{x}, \bm{y} \in \mathring{\bm{M}}$, $\bm{Z}(\bm{x},\bm{y}) = \tilde{Z}(x,y)= \overline{Z}(x,y,z)$.
    \begin{align}
        0 \ge  &4 \frac{\varphi''}{\bm{L}} + ( 1- \varphi'^{2} ) \int_{\alpha} K + 2 \Big( \varphi - \varphi' \frac{\bm{l}}{\bm{L}} \Big)  \int_{\Gamma_{t}} \kappa^{2} ds + \varphi' \int_{[\bm{\Gamma}_{t}(x):\bm{\Gamma}_{t}(y)]} \kappa^{2} ds - \bm{L} \partial_{t} \varphi \\ \nonumber &+ \frac{2\kappa^{\partial N}(z)} {\Big\langle \frac{\alpha'_{+}(t_{0})-\alpha'_{-}(t_{0})}{\tilde{d}},\mathcal{N}^{\partial N}(z) \Big\rangle}.
    \end{align}
    
\end{enumerate}
\end{prop}

We now discuss the lower bounds of the chord-arc profile on surfaces. The following lemma gives rise to the estimate of total curvature between two points along the boundary when the distance between two points is small. For $z_{1}, z_{2} \in \partial N$, denote $[z_{1}:z_{2}]$ as the smaller portion of $\partial N \setminus \{ z_{1}, z_{2} \}$. 
\begin{lem} There exists $\epsilon_{0} = \epsilon_{0}(N)>0$ such that the following holds: Suppose $\epsilon \in (0, \epsilon_{0})$. Let $\Gamma$ be a curve in $(N,\partial N,g)$ which meets $\partial N$ orthogonally at $\partial \Gamma = \{ z_{0}, z_{1} \}$ with length $L$. Denote $C = \sup_{\partial N \cap B_{3L}(\partial \Gamma)} \kappa^{\partial N}$. If $L(1+C) \le \min (\frac{\epsilon}{100}, \frac{|\partial N|}{8})$ and $z \in \partial N$ is a point achieving $\tilde{d}(x,y) = d(x,z)+d(z,y)$ for some $x,y \in \Gamma$, then 
\begin{equation}
    \int_{[z_{0}:z_{1}]} \kappa \le \frac{2\epsilon}{5} \text{ and } \int_{[z_{0}:z]} \kappa \le \frac{2\epsilon}{5}.
\end{equation}
\begin{proof} 
We regard $N$ as a convex domain of a closed manifold $\hat{N}$. We parametrize $\Gamma$ with arclength by $\Gamma : [0,|\Gamma|] \rightarrow N$. Fix a point $z_{0} \in \partial N \cap \Gamma$ and consider an exponential map $\exp_{z_{0}}: B_{a}(0) \rightarrow \hat{N}$ for some $a>0$. We also consider a geodesic normal coordinate $x = (x_{1},x_{2})$ such that $x  := X \circ \exp_{z_{0}}^{-1}$ where $X = (X_{1},X_{2})$ is an Euclidean coordinate. In the geodesic normal coordinate, note that $g_{ij} = \delta_{ij} + O(r^{2})$ and Christoffel symbols satisfy $\Gamma_{ij}^{k} = O(r)$ where $r = \sqrt{X_{1}^{2}+X_{2}^{2}}$. Denote $\theta(s) = \tan^{-1}(\dot{x}_{2}(\Gamma(s))/\dot{x}_{1}(\Gamma(s)))$ as an angle in the geodesic normal coordinate. We consider the local geodesic equation $|\nabla_{\dot{\Gamma}}\dot{\Gamma}|^{\perp} = \kappa_{\Gamma} N_{\Gamma}$ in terms of local coordinates:
\begin{equation}
    \sum_{k=1,2} \Big( \ddot{x}_{k}+\sum_{i,j=1,2}\dot{x}_{i}\dot{x}_{j}\Gamma_{ij}^{k} \Big)\frac{\partial}{\partial x_{k}} =  (\kappa+O(r)) \Big(-\dot{x}_{2} \frac{\partial}{\partial x_{1}}+\dot{x}_{1} \frac{\partial}{\partial x_{2}}\Big).
\end{equation}
    From (34), we obtain
    \begin{equation}
        \frac{\partial \theta(s)}{\partial s} = \frac{\ddot{x}_{2} \dot{x}_{1} - \ddot{x}_{1} \dot{x}_{2}}{\dot{x}^{2}_{1}+\dot{x}^{2}_{2}} = \kappa(s) + O(r).
    \end{equation}
By (35) and the definition of the geodesic normal coordinate, there exists $\epsilon_{0} = \epsilon_{0}(N)>0$ such that the arguments in Lemma 3.5 in \cite{ACGL} to bound the length of $[z_{0}:z_{1}]$ in terms of $d(z_{0},z_{1})$ works directly. We follow the arguments in the proof of Lemma 5.2 in \cite{LZ} and have
\begin{equation}
    |[z_{0}:z_{1}]| \le \frac{2\epsilon}{5C} \text{ and }
    |[z_{0}:z]| \le \frac{2\epsilon}{5C}.
\end{equation}
Note that the constant in the proof of Lemma 5.2 in \cite{LZ} is not sharp and we took a more strict upper bound. Then we finally obtain
\begin{equation*}
    \int_{[z_{0}:z_{1}]} \kappa^{\partial N} ds \le C |[z_{0}:z_{1}]| \le C \cdot \frac{2\epsilon}{5C} = \frac{2\epsilon}{5}
\end{equation*}
and another inequality of (33) also follows.
\end{proof}
\end{lem}
Note that there exists $L_{0} = L_{0}(N,g)$ such that for a topological disk $A \subseteq N$ with $|\partial A| \le L_{0}$, the isoperimetric inequality holds $C' |A| \le |\partial A|^{2}$ for some explicit constant $C'$ by Proposition 2.1 in \cite{BR}. Denote $K_{0} = \sup_{x \in N} |K|$.
\begin{thm}
    Suppose $L(t) \rightarrow 0$ as $t \rightarrow T$. Take $\epsilon_{1}= \epsilon_{1}(N)>0$. Given any $\epsilon \in (0,\epsilon_{1})$, there exists $c_{\epsilon}$ such that the following holds: let $\{ \Gamma_{t} \}_{t \in [0,T)}$ be a free boundary curve shortening flow on $N$. Suppose $L(0) (1+C) \le \min(\frac{\epsilon}{100},\frac{L_{0}}{3})$, where $C = \sup_{\partial N} \kappa^{\partial N}$. Given any $c_{0} \in (0,c_{\epsilon})$, if the inequality
    \begin{equation}
        \bm{\psi}(\delta,t) \ge \begin{cases}
            c_{0} \bm{L}(t) \Big\{ \sin \Big((\pi-\epsilon) \frac{\delta}{\bm{L}(t)} + \frac{\epsilon}{2} \Big)+ 64 (\frac{\delta}{\bm{L}(t)} -\frac{1}{4})^{3} \sin \frac{\epsilon}{2} \Big\} e^{-K_{0}t} &\text{ if }0 \le \delta \le \frac{\bm{L}(t)}{4} \\
              c_{0} \bm{L}(t) \sin \Big((\pi-\epsilon) \frac{\delta}{\bm{L}(t)} + \frac{\epsilon}{2} \Big) e^{-K_{0}t} &\text{ otherwise}.
        \end{cases}
    \end{equation}
    holds at $t=0$, then it holds for all $t \in [0,T)$.
    \begin{proof}
        We take $\epsilon_{1} = \min(\epsilon_{0},L_{0},\frac{2}{C'K_{0}},\frac{\pi}{20})$. Take $\varphi \in C^{2}([0,\frac{1}{2}] \times [0,T))$ as following:
        \begin{equation}
            \varphi(\zeta,t) = \begin{cases}  c_{0} \{ \sin ((\pi-\epsilon) \zeta + \frac{\epsilon}{2} )+ 64 (\zeta -\frac{1}{4})^{3} \sin \frac{\epsilon}{2} \} e^{-K_{0}t} &\text{ if }\zeta \in [0,\frac{1}{4}]\\
             c_{0} \{ \sin ((\pi-\epsilon) \zeta + \frac{\epsilon}{2} )\} e^{-K_{0}t} &\text{ otherwise}
            \end{cases}
        \end{equation}
        and extend (38) to $\varphi \in C^{2}([0,1] \times [0,T))$ to satisfy $\varphi(1-\zeta,t) = \varphi(\zeta,t)$ for $\zeta \in [0,1/2]$ and $t \in [0,T)$. We will choose $c_{0}>0$ later in the proof. Note that $\varphi(0,t) = 0$ and $\partial_{\zeta} \varphi (1/2,t) =0$ for every $t \in [0,T)$. We define the auxiliary function $\bm{Z}$ with $\varphi$ defined in (38). Notice that the initial condition holds by (37) and we argue by contradiction. Denote $t_{0}:= \sup \{ t \in [0,T): \bm{Z}(\cdot,\cdot,t) \ge 0 \}$ and assume $t_{0} <T$. Suppose $\bm{Z}(\bm{x},\bm{y},t_{0})=0$ and note that this satisfies the conditions of Proposition 4.1. Denote $\Gamma_{t_{0}} \cap \partial N = \{ z_{0}, z_{1} \}$. First we estimate
    \begin{equation*}
       \Theta :=\int_{\Gamma_{t_{0}}} \kappa ds \text{  and  } \omega:=\int_{[\bm{\Gamma}_{t_{0}}(x):\bm{\Gamma}_{t_{0}}(y)]} \kappa ds.
    \end{equation*}
    Denote the region surrounded by $[z_{0}:z_{1}]$ and $\Gamma_{t_{0}}$ by $A_{t_{0}}$.
   
    Since $\Gamma_{t_{0}}$ meets $\partial N$ orthogonally, by Gauss-Bonnet theorem we have
    \begin{align}
        \nonumber \Theta &= 2 \pi - \int_{[z_{0}:z_{1}]} \kappa^{\partial N} - \int_{A_{t_{0}}} K - 2 \cdot \frac{\pi}{2} \\
        &\ge \pi - \frac{2\epsilon}{5} - \int_{A_{t_{0}}} K ,
    \end{align}
    where we obtain (39) by applying Lemma 4.2.

    Now we deduce the lower bound of $\omega$ and we consider the case of $\text{sign}(\bm{x}) = \text{sign}(\bm{y})$ first. By Proposition 3.1, $A_{(\bm{x},\bm{y})}$ is a topological disk. We denote $\beta = \cos ^{-1} \varphi'$ and $\beta$ is an interior angle between two curves at $x$ and $y$ of $A_{(\bm{x},\bm{y})}$ by (7). We obtain the following by Gauss-Bonnet theorem:
    \begin{align}
    \nonumber \omega=\int_{[\bm{\Gamma}_{t_{0}}(\bm{x}):\bm{\Gamma}_{t_{0}}(\bm{y})]} \kappa ds &= 2 \pi - 2(\pi-\beta) - \int_{A_{(\bm{x},\bm{y})}} K \\ &= 2 \cos ^{-1} \varphi' - \int_{A_{(\bm{x},\bm{y})}} K.
    \end{align}

    Now we consider the case $\text{sign}(\bm{x}) \neq \text{sign}(\bm{y})$. We apply Proposition 3.1 again and obtain that $A_{(\bm{x},\bm{y})}$ is a topological disk in $\tilde{N}$. Moreover, denote $\beta = \cos ^{-1} \varphi'$ to be an interior angle between two curves $[\bm{\Gamma}_{t_{0}}(\bm{x}):\bm{\Gamma}_{t_{0}}(\bm{y})]$ and $\bm{\alpha}$ at $\bm{x}$ and $\bm{y}$ by (16). We separate $A_{\bm{\Gamma}_{t_{0}}}$ by
    \begin{equation*}
        A_{1} = A_{(\bm{x},\bm{y})} \cap N \text{ and } A_{2} = A_{(\bm{x},\bm{y})} \cap (\tilde{N} \setminus N).
    \end{equation*}Note that $A_{1}$ and $A_{2}$ are both topological disks. Without loss of generality, $z_{0} \in [\Gamma_{t_{0}}(x):\Gamma_{t_{0}}(y)]$. Since $\Gamma_{t_{0}}$ orthogonally meets $\partial N$, by applying (18) and Gauss-Bonnet theorem we have
    \begin{align}
    \nonumber \omega=\int_{[\bm{\Gamma}_{t_{0}}(x):\bm{\Gamma}_{t_{0}}(y)]} \kappa ds &= \int_{[\bm{\Gamma}_{t_{0}}(x):z_{0}]} \kappa ds + \int_{[z_{0}:\bm{\Gamma}_{t_{0}}(y)]} \kappa ds \\ \nonumber &= \Big( 2 \pi - \theta_{0} - \pi/2 - (\pi-\beta) - \int_{A_{1}} K  - \int_{[z_{0}:z]} \kappa^{\partial N} \Big) \\ \nonumber &+ \Big( 2 \pi - (\pi-\theta_{0}) - \pi/2 - (\pi-\beta) - \int_{A_{2}} K  - \int_{[z_{0}:z]} \kappa^{\partial N} \Big) \\ \nonumber &= 2 \cos ^{-1} \varphi' - \int_{A_{(\bm{x},\bm{y})}} K - 2 \int_{[z_{0}:z]} \kappa^{\partial N} \\ &\ge 2 \cos ^{-1} \varphi' - \int_{A_{(\bm{x},\bm{y})}} K - \frac{4\epsilon}{5},
    \end{align}
    where we obtain (41) by applying Lemma 4.2.
     
    By the isoperimetric inequality and our choice of $L_{0}$, we obtain the following estimate of the area of $A_{t_{0}}$ and $A_{(\bm{x},\bm{y})}$:
    \begin{align}
        \nonumber |A_{t_{0}}| &\le  C' ( L(t_{0}) +|[z_{0}:z_{1}]| )^{2} \\ &\le C' (\frac{\epsilon}{100} + \frac{\epsilon}{2})^{2} \le C' \epsilon^{2}, \\
        |A_{(\bm{x},\bm{y})}| &\le C' ( 2L(t_{0}))^{2} \le \frac{C' \epsilon^{2}}{100}.
    \end{align} We obtain (42) by applying (36) and (43) by applying $|\alpha| \le [\bm{\Gamma}_{t_{0}}(x):\bm{\Gamma}_{t_{0}}(y)] \le L(t_{0})$. Moreover, by (42) and (43), we have the control of total Gaussian curvature of $A_{t_{0}}$ and $A_{(\bm{x},\bm{y})}$:
    \begin{align}
        \int_{A_{t_{0}}} K &\le K_{0} |A_{t_{0}}| \le C'K_{0}\epsilon^{2}, \\ \int_{A_{(\bm{x},\bm{y})}} K &\le  K_{0} |A_{(\bm{x},\bm{y})}| \le \frac{1}{100}C'K_{0}\epsilon^{2}.
    \end{align}
    By our choice of $\epsilon_{1}$, (39)-(41), (44) and (45), we have
    \begin{equation}
        \Theta \ge \pi -\epsilon \text{ and } \omega \ge 2 \Big(\cos ^{-1} \varphi'- \frac{\epsilon}{2}\Big).
    \end{equation}
    By our choice of $\varphi$, $|\varphi'| \le c_{0} \{ (\pi- \epsilon) + 12 \sin \frac{\epsilon}{2} \}$ holds. By taking sufficiently small $c_{0}$, we have $\cos ^{-1} \varphi' \ge \frac{\epsilon}{2}$ and our estimate of $\omega$ in (46) is proper.
    We apply Hölder inequality and obtain
    \begin{align}
        \int_{\Gamma_{t}} \kappa^{2} ds &\ge |\Gamma_{t}|^{-1} \Big( \int_{\Gamma_{t}} |\kappa| ds \Big)^{2} \ge \frac{2}{\bm{L}} \Theta^{2} \ge \frac{2}{\bm{L}} (\pi-\epsilon)^{2}
    \end{align}
    and
    \begin{align}
       \nonumber \int_{[\bm{\Gamma}_{t}(x):\bm{\Gamma}_{t}(y)]} \kappa^{2} ds &\ge |[\bm{\Gamma}_{t}(x):\bm{\Gamma}_{t}(y)]|^{-1} \Big( \int_{[\bm{\Gamma}_{t}(x):\bm{\Gamma}_{t}(y)]} |\kappa| ds \Big)^{2} \\ &\ge \frac{1}{l} \omega^{2} \ge \frac{4}{l} \Big( \cos ^{-1} \varphi'- \frac{\epsilon}{2} \Big)^{2}.
    \end{align}
    Note that We apply (47) and (48) to (31) and (32), in either case we have
    \begin{align}
 \nonumber 0 &\ge  4 \frac{\varphi''}{\bm{L}} + ( 1- \varphi'^{2} ) \int_{\alpha} K + 2 \Big( \varphi - \varphi' \frac{\bm{l}}{\bm{L}} \Big)  \int_{\Gamma_{t}} \kappa^{2} ds + \varphi' \int_{[\Gamma_{t}(x):\Gamma_{t}(y)]} \kappa^{2} ds - \bm{L} \partial_{t} \varphi \\ \nonumber &\ge 4 \frac{\varphi''}{\bm{L}} + ( 1- \varphi'^{2} ) \int_{\alpha} K + \frac{4}{\bm{L}} \Big( \varphi - \varphi' \frac{\bm{l}}{\bm{L}} \Big)  ( \pi-\epsilon)^{2} + \frac{4}{l} \Big( \cos ^{-1} \varphi'- \frac{\epsilon}{2} \Big)^{2} \varphi' - \bm{L} \partial_{t} \varphi \\ \nonumber
  &\ge 4 \frac{\varphi''}{\bm{L}} - K_{0} d + \frac{4}{\bm{L}} \Big( \varphi - \varphi' \frac{\bm{l}}{\bm{L}} \Big)  (\pi-\epsilon)^{2} + \frac{4}{l} \Big( \cos ^{-1} \varphi'- \frac{\epsilon}{2} \Big)^{2}  \varphi' - \bm{L} \partial_{t} \varphi \\&= 4 \frac{\varphi''}{\bm{L}} - K_{0} \bm{L} \varphi + \frac{4}{\bm{L}} \Big( \varphi - \varphi' \frac{\bm{l}}{\bm{L}} \Big)  (\pi-\epsilon)^{2} + \frac{4}{l} \Big( \cos ^{-1} \varphi'- \frac{\epsilon}{2} \Big)^{2} \varphi' - \bm{L} \partial_{t} \varphi.
    \end{align}
If $\epsilon \in (0,\pi/20)$, then the following holds by direct calculation:
\begin{equation}
    (\pi-\epsilon) \cos \Big( (\pi-\epsilon)\frac{1}{4} + \frac{\epsilon}{2} \Big) > 12 \sin \frac{\epsilon}{2}.
\end{equation}
Moreover, since $\cos ^{-1} \varphi' \rightarrow \pi/2$ uniformly as $c_{0} \rightarrow 0$ in $[0,1/4]$, we can take $c_{0}$ sufficiently small to satisfy
\begin{equation}
    \Big( \cos ^{-1} \varphi'- \frac{\epsilon}{2} \Big)^{2}-(\pi-\epsilon)^{2} \zeta^{2} \ge 1
\end{equation}
for $\zeta \in [0,1/4]$.
Let us consider the following term for $\zeta \in (0,1/4]$:
\begin{align}
    \nonumber &\varphi''+ (\pi-\epsilon)^{2} \varphi + \zeta^{-1} \varphi'\Big(  \Big( \cos ^{-1} \varphi'- \frac{\epsilon}{2} \Big)^{2}-(\pi-\epsilon)^{2} \zeta^{2} \Big) \\ \nonumber
    &= c_{0} e^{-K_{0}t}\Big( 384 \Big(\zeta-\frac{1}{4}\Big)\sin \frac{\epsilon}{2} \Big) + c_{0} \zeta^{-1} e^{-K_{0}t} \Big( (\pi-\epsilon) \cos \Big( (\pi-\epsilon)\zeta + \frac{\epsilon}{2} \Big) + 192 \Big(\zeta-\frac{1}{4}\Big)^{2}\sin \frac{\epsilon}{2}\Big)\\ \nonumber&\Big(  \Big( \cos ^{-1} \varphi'- \frac{\epsilon}{2} \Big)^{2}-(\pi-\epsilon)^{2} \zeta^{2} \Big) \\ &> c_{0} e^{-K_{0}t} \Big( 384 \Big(\zeta-\frac{1}{4}\Big)\sin \frac{\epsilon}{2} \Big) + c_{0} \zeta^{-1} e^{-K_{0}t} \Big( 12 \sin \frac{\epsilon}{2} + 192 \Big(\zeta-\frac{1}{4}\Big)^{2}\sin \frac{\epsilon}{2}\Big) \\
    \nonumber &\ge c_{0} e^{-K_{0}t} \Big( 384 \Big(\zeta-\frac{1}{4}\Big)\sin \frac{\epsilon}{2} \Big) + 4c_{0}  e^{-K_{0}t} \Big( 12 \sin \frac{\epsilon}{2} + 192 \Big(\zeta-\frac{1}{4}\Big)^{2}\sin \frac{\epsilon}{2}\Big) \\
    &= 48 c_{0} e^{-K_{0}t} \sin \frac{\epsilon}{2} \Big( 4 \Big(\zeta-\frac{1}{4}\Big) + 1\Big)^{2} \ge 0.
\end{align}
(52) comes from (50), (51) and the monotone decreasing property of the Cosine function. For $\zeta \in [1/4,1/2]$, both of the following holds for small $c_{0}>0$:
\begin{align}
    \varphi''+ (\pi-\epsilon)^{2} \varphi &=0 \\
    \Big( \cos ^{-1} \varphi'- \frac{\epsilon}{2} \Big)^{2}-(\pi-\epsilon)^{2} \zeta^{2} &>0.
\end{align}

We also have the following identity in time derivative terms:
\begin{equation}
    K_{0} \partial_{t} \varphi  + \varphi =0.
\end{equation}
By applying (53)-(56) into (49), we obtain the contradiction.
    \end{proof}
\end{thm}
\begin{thm}
Suppose $L(t) \nrightarrow 0$ as $t \rightarrow T$ where $T < \infty$. Let $\{ \Gamma_{t} \}_{t \in [0,T)}$ be a free boundary curve shortening flow on $N$. If $\bm{L}_{T} := \lim_{t \rightarrow T} \bm{L}(t)>0$, then if the inequality
\begin{equation*}
    \bm{\psi}(\delta,t) \ge c_{0} \bm{L}(t) e^{\big(-\frac{4 \pi^{2}}{L_{T}^{2}} - K_{0}\big)t} \sin \Big(\frac{\pi \delta}{\bm{L}(t)} \Big)
\end{equation*}
holds at $t=0$, then it holds for all $t \in [0,T)$.
\begin{proof}
We modify the proof of Theorem 5.4 in \cite{LZ}. We adopt the modified time coordinate $\tau := \int^{t}_{0} \frac{1}{\bm{L}(s)^{2}} ds$. We take $\varphi: [0,1] \times [0,T)$ by
\begin{equation*}
    \varphi(\zeta,t) = c_{0} e^{-4\pi^{2}\tau(t)-K_{0}t}\sin (\pi \zeta).
\end{equation*}
As before, denote $t_{0}:= \sup \{ t \in [0,T): \bm{Z}(\cdot,\cdot,t) \ge 0 \}$ and assume $t_{0}<T$. Then there exists $(\bm{x} , \bm{y}) \in (\bm{M} \times \bm{M} ) \setminus \bm{D}$ such that $\bm{Z}(\bm{x}_{0},\bm{y}_{0},t_{0}) = \min_{(\bm{x},\bm{y}) \in (\bm{M} \times \bm{M} ) \setminus \bm{D}} \bm{Z}(\bm{x},\bm{y},t_{0})$.  By Proposition 4.1, in either case we have
\begin{align}
      \nonumber 0 &\ge  4 \frac{\varphi''}{\bm{L}} + ( 1- \varphi'^{2} ) \int_{\alpha} K + 2 \Big( \varphi - \varphi' \frac{\bm{l}}{\bm{L}} \Big)  \int_{\Gamma_{t}} \kappa^{2} ds + \varphi' \int_{[\Gamma_{t}(x):\Gamma_{t}(y)]} \kappa^{2} ds - \bm{L} \partial_{t} \varphi \\ \nonumber &\ge  4 \frac{\varphi''}{\bm{L}} - K_{0}\bm{L}\varphi + 2 \Big( \varphi - \varphi' \frac{\bm{l}}{\bm{L}} \Big)  \int_{\Gamma_{t}} \kappa^{2} ds + \varphi' \int_{[\Gamma_{t}(x):\Gamma_{t}(y)]} \kappa^{2} ds - \bm{L} \partial_{t} \varphi \\ &>  4 \frac{\varphi''}{\bm{L}} - K_{0}\bm{L}\varphi  - \bm{L} \partial_{t} \varphi,
\end{align}
where (57) follows from the strict concavity of $\varphi(\cdot,t)$. But our choice of $\varphi$ gives
\begin{align}
\nonumber 4 \frac{\varphi''}{\bm{L}} - K_{0}\bm{L}\varphi  - \bm{L} \partial_{t} \varphi &= -4 \pi^{2} \frac{\varphi}{\bm{L}} - K_{0}\bm{L}\varphi + \bm{L}(4 \pi^{2} \tau'(t) +K_{0}) \varphi \\ &= -4 \pi^{2} \frac{\varphi}{\bm{L}} - K_{0}\bm{L}\varphi + \bm{L}\Big(4 \pi^{2} \frac{1}{\bm{L}^{2}} +K_{0}\Big) \varphi =0
\end{align}
(57) and (58) give the contradiction. And since $\bm{L}$ is a decreasing function in time $t$, the claim of the Theorem follows. 
\end{proof}
\end{thm}
Let us denote $\lambda : M \rightarrow \mathbb{R}$ by the arclength to the nearer endpoint in the sense of internal distance in $\Gamma_{t}$. Together with the proof of Proposition 5.5 in \cite{LZ}, Theorem 4.3 and Theorem 4.4 gives the following boundary avoidance estimate.
\begin{prop} 
Let $\{ \Gamma_{t} \}_{t \in [0,T)}$ be a free boundary curve shortening flow on $N$. Given any $\delta>0$, there exists $\epsilon = \epsilon(\Gamma_{0},N,\delta)>0$ such that
\begin{equation*}
    \lambda(x,t) > \delta \Rightarrow d(\gamma(x,t), \partial N)> \epsilon.
\end{equation*}
\end{prop}
As in \cite{LZ}, Proposition 4.5 and arguments of the proof of Theorem 6.1 in \cite{LZ} gives the following long-time behavior.
\begin{thm} Let $(N,\partial N, g)$ be a closed Riemannian surface with convex boundary and $\{ \Gamma_{t} \}_{t \in [0,T)}$ be a maximal free boundary curve shortening flow starting from a properly embedded closed interval $\Gamma_{0}$ in $N$. Then either:
\begin{enumerate}
    \item $T = \infty$, in which case $\Gamma_{t}$ converges smoothly as $t \rightarrow \infty$ to an embedded geodesic in $N$ which meets $\partial N$ orthogonally; or
    \item $T< \infty$, in which case $\Gamma_{t}$ converges uniformly to some single half-round point $z \in \partial N$ smoothly in the sense of the blow up limit of the curve converges to the unit semi-circle.
\end{enumerate}
\end{thm}
\section{The family of curves and tightening procedure}
In this section, we formulate the min-max construction of free boundary embedded geodesics on Riemannian $2$-disks with convex boundary. We discuss the smooth min-max setting and obtain the existence result via proper pull-tight procedure. Theorem B in \cite{Sm} (see also Appendix of Hatcher's work \cite{Ha}) proves that the space of (unparametrized) embedded intervals on $D^{2}$ whose endpoints are on $\partial D^{2}$ relative to the space of point curves retracts onto $\mathbb{R}P^{2}$. Let us denote the space of embedded curves by $\Sigma$ and denote by $\Sigma_{0}$ the space of point curves. We denote $\mathcal{S} = \Sigma/\Sigma_{0}$. We consider two distinct nontrivial relative homology classes $h_{1}$ and $h_{2}$ on the space of embedded intervals:
\begin{equation*}
    h_{i} := H_{i}(\mathcal{S}, \mathbb{Z}_{2}) = \mathbb{Z}_{2}.
\end{equation*}
Let us consider the $\mathbb{Z}_{2}$-Cohomology ring $H^{*}(\mathcal{S},\mathbb{Z}_{2})$ of $\mathcal{S}$. We denote $\alpha$ to be a generator of the first cohomology ring. Then the cohomology ring is:
\begin{equation*}
    H^{*}(\mathcal{S},\mathbb{Z}_{2}) = \mathbb{Z}_{2}[\alpha]/\alpha^{2}.
\end{equation*}

Let $IV_{1}(D^{2})$ be a space of integral varifolds on $(D^{2},g)$ and endow an $F$-metric on the space of varifolds as in 2.1(19)(20) in \cite{P}.

For each $i=1,2$, we now define the $i$-sweepout. We denote an $i$-dimensional simplicial complex by $X$. If
\begin{equation*}
    \Phi^{*} (\omega) \neq 0,
\end{equation*}
then we say that $\Phi:X \rightarrow \mathcal{S}$ detects $\omega \in H^{i}(\mathcal{S}, \mathbb{Z}_{2})$. We define $\Phi$ to be an \emph{$i$-sweepout} endowed with a smooth topology if $\Phi$ detects $i$-th cup product $\alpha^{i} \in H^{i}(\mathcal{S},\mathbb{Z}_{2})$. We define the \emph{width} of $i$-parameter sweepouts as
\begin{equation}
    \omega_{i}(D^{2}) : = \inf _{\Phi \in S_{i}} \sup_{x \in X} |\Phi(x)| = L_{i},
\end{equation}
for $i \in \{ 1, 2 \}$. By the definition of $i$-sweepouts in (59), $\omega_{1}(D^{2}) \le \omega_{2}(D^{2})$ holds. 

We define a \emph{minimizing sequence} to be the sequence of $i$-sweepouts such that $\lim_{j \rightarrow \infty} \sup_{x \in X} |\Phi_{j}(x)| = L_{i}$. We denote sequence of curves $\Phi_{j}(x_{j})$ to be a min-max sequence if $|\Phi_{j}(x_{j})|$ converges to $L_{i}$, where $x_{j} \in X$ and $\{ \Phi_{j}(x) \}$ is a minimizing sequence. We define the \emph{critical set} $\Lambda(\{ \Phi_{j} \})$ to be a set of stationary varifolds achieved by the limit of min-max sequence induced by  $\{ \Phi_{j}(x) \}$. We denote the set of critial geodesic by $W_{L_{i}}$ which is  a set of stationary varifolds whose support is a free boundary embedded geodesic and length is $L_{i}$.

We follow Abraham's proof of the bumpy metric theorem for curves \cite{Ab} and Theorem 9 in Ambrozio-Carlotto-Sharp \cite{ACS} of the free boundary minimal surface version which proves the genericity of bumpy metric in $C^{r}$-Baire sense. The compactness of free boundary embedded geodesics with bounded length follows from arguments in Appendix A in \cite{K2}.

By tightening argument via free boundary curve shortening flow of Theorem 4.6 on long-time behavior of free boundary curve shortening flow, and applying the arguments in the proof of Theorem 2.1 in \cite{K1}, we have the following existence of free boundary embedded geodesics achieving the width:
\begin{thm}
Suppose $(D^{2},\partial D^{2},g)$ to be a smooth Riemannian $2$-disk with a strictly convex boundary. For $i=1,2$ and any minimizing sequence $\{ \Phi_{j} \}$ of $i$-sweepouts, there is a deformed minimizing sequence $\{\hat{ \Phi}_{j} \}$ of $\{ \Phi_{j} \}$ satisfying the following property. For any $s>0$, there is some $0<a<L_{i}$ satisfying
\begin{equation*}
        \{ \hat{\Phi}_{j}(x) \in IV_{1}(D^{2}) : |\hat{\Phi}_{j}(x)| \ge L_{i}-a \} \subset \bigcup_{\gamma \in \Lambda(\{ \Phi_{j} \})\cap W_{L_{i}}} B^{F}_{s}(\gamma)
        \end{equation*}
for all sufficiently large $j$, where $B^{F}_{s}(\gamma)$ is a $F$-metric ball with center $\gamma$. Moreover, the multiplicity of geodesics in the critical set is $1$.
\end{thm}
\begin{rmk}
We still can run the tightening procedure even initial curves do not meet orthogonally with the boundary $\partial D^{2}$. Geometrically, we can slightly deform the curves near the boundary to make the intersection angles the right angle. Moreover, we can apply the flow in the `weak' sense. We may consider the free boundary curve shortening flow as a Neumann boundary problems for nonlinear parabolic PDEs (See Chapter 10 of \cite{Fr}). Then the weak solution of the flow has a orthogonal boundary condition and smoothness at any positive time before the maximum existence time $t \in (0,T)$. The flow is still a $C^{0}$-solution at $t \in [0,T)$. 
\end{rmk}

By applying the classical Lusternik-Schnirelmann argument, we prove that if two widths are the same, then there exists an $S^{1}$-cycle of free boundary embedded geodesics. This gives the existence of two free boundary geodesics on any Riemannian disk with a strictly convex boundary. Moreover, if the metric on $(D^{2},\partial D^{2},g)$ is bumpy, then we deduce the existence of two free boundary embedded geodesics with distinct lengths. The argument in the proof of Corollary 2.2 in \cite{K1} applies directly.
\begin{thm}
Suppose $(D^{2},\partial D^{2},g)$ to be a Riemannian $2$-disk with strictly convex boundary. If $w_{1}(D^{2}) = w_{2}(D^{2})$, then there exist infinitely many distinct free boundary embedded geodesics in $(D^{2},\partial D^{2}, g)$.
\end{thm}
\begin{cor}
 On the Riemannian disk endowed with a bumpy metric and strictly convex boundary, there are at least two free boundary embedded geodesics with length $L_{1}$ and $L_{2}$.
\end{cor}
\section{Morse Index Bound}
In this section, we obtain the generic Morse Index bound of free boundary embedded geodesics on Riemannian $2$-disk. We mainly follow the idea in Section 6 and Section 7 of \cite{K1} which proves the Morse Index of simple closed geodesics on bumpy spheres based on the interpolation technique based on quantitative $F$-distance estimate. We focus on the necessary modification to prove the Morse Index bound in the free boundary setting. Throughout this section, we suppose that $(D^{2},\partial D^{2},g)$ is endowed with a bumpy metric. Suppose $D^{2}$ is embedded in some closed surface $\tilde{D}^{2}$.

Let us fix a free boundary embedded geodesic $\gamma$. We adopt the Fermi coordinate $c: [0,L] \times (-h,h) \rightarrow D^{2}$ on $N_{h}(\gamma)$ on the tubular neighborhood $N_{h}(\gamma)$ of the fixed geodesic $\gamma$. Moreover, we adopt the metric perturbation in Section 5.2 in \cite{K1} and follow the notation therein. By choosing a sufficiently small $\beta>0$ and sufficiently large $M > \sup_{D^{2}} |K|$ in Proposition 5.2 of \cite{K1} and applying the corresponding deformation (17) in \cite{K1}, we obtain the strict stability of the free boundary embedded geodesic $\gamma_{g_{\beta}}$ in the perturbed metric.

\begin{prop} There exists small $\beta>0$ and $M>\sup_{D^{2}}|K|$ satisfying the following: $\gamma_{g_{\beta}}$ is a strictly stable geodesic and the ambient Gaussian curvature is strictly negative.
\begin{proof}
Denote $\partial \gamma_{g_{\beta}} = \{p_{1},p_{2} \}$. Then by the change of the second fundamental form by conformal deformation following \cite{Be}, for $i=1,2$, we have a geodesic curvature of boundary on $\kappa^{\partial D^{2}}_{g_{\beta}}$ as
\begin{equation*}
    \kappa^{\partial D^{2}}_{g_{\beta}}(p_{i}) = e^{-\phi_{\beta}} \Big( \kappa^{\partial D^{2}}_{g}(p_{i}) - \frac{\partial \phi_{\beta}}{\partial \nu }\Big) = e^{-\phi_{\beta}} \kappa^{\partial D^{2}}_{g}(p_{i}).
\end{equation*}
Also note that $K_{g_{\beta}}(x) = K(x)-M<0$ on $x \in \gamma_{g_{\beta}}$ by (19) of \cite{K1} and this proves the latter conclusion of the claim. Now we prove that we can choose $\beta$ and $M$ such that the second variation (3) of $\gamma_{g_{\beta}}$ is positive definite. It suffices to show that by taking the suitable $\beta$ and $M$, the eigenvalue $\{ \lambda_{k} \}$ with associated eigenfunctions $\{ \phi_{k}\}$ of the following equation with Robin boundary condition from (4) are all positive:
\begin{equation}
    \begin{cases} (\Delta_{\gamma_{g_{\beta}}} + K-M) \phi_{k} + \lambda_{k} \phi_{k} = 0 \text{ on }\gamma_{g_{\beta}}\\
    \frac{\partial \phi_{k}}{\partial \eta}(p_{i}) -e^{-\phi_{\beta}}\kappa (p_{i}) \phi_{k} (p_{i}) = 0 \text{ for } i=1,2.
    \end{cases}
\end{equation}
By the condition (iii) in Proposition 5.2 in \cite{K1}, we can take $e^{-\phi_{\beta}}$ to be arbitrarily close to $1$ by taking sufficiently small $\beta$. By taking sufficiently large $M$ and sufficiently small $\beta$, we can obtain all the eigenvalues $\{ \lambda_{k} \}$ to be strictly positive by standard elliptic theory. 
\end{proof}
\end{prop}

We now discuss the free boundary mean convex neighborhood of the geodesic $\gamma$ which is strictly stable and with negative ambient Gaussian curvature to adopt the squeezing lemma in \cite{K1}. We adopt the idea in Proposition 2.4 of \cite{HK} which constructs the free boundary mean convex neighborhood of the free boundary embedded minimal surface in $3$-ball with a strict convex boundary. Consider the first eigenvalue $\lambda_{1} >0$ and the associated eigenfunction $\phi_{1} \in C^{\infty}(\gamma)$ with $\int_{\gamma} \phi_{1}^{2} ds=1$ which is a solution of the equation (4). Without loss of generality, we can set $\phi_{1}$ to be strictly positive in $\gamma$. We define
\begin{equation}
C_{\gamma} := \frac{\max_{\gamma} \phi_{1}}{\min_{\gamma} \phi_{1}} \ge 1
\end{equation}
and call $C_{\gamma}$ by \emph{Harnack constant} of $\gamma$.
\begin{prop}
    Suppose $\gamma$ is a strictly stable geodesic on a Riemannian $2$-disk $(D^{2},\partial D^{2},g)$ with strictly convex boundary and with negative ambient Gaussian curvature. Then there is a neighborhood $N_{h}(\gamma)$ foliated by a free boundary mean convex foliation $\{\gamma_{t} \}_{t \in [-\epsilon,\epsilon]}$ of $\gamma$ satisfying the following:
    \begin{enumerate}
        \item $\gamma_{0} = \gamma$,
        \item $\gamma_{t}$ has a mean curvature vector towards $\gamma_{0}$ for $t \in [-\epsilon,\epsilon]$,
        \item $C_{t}:= \max{d_{x \in \gamma_{t}}(x,\gamma)}/\min{d_{x \in \gamma_{t}}(x,\gamma)} \le 2C_{\gamma}$ for $t \in [-\epsilon,0) \cup (0,\epsilon]$.
    \end{enumerate}
\begin{proof}
Let us denote the unit normal vector on $\gamma$ by $\nu$. Then we consider normal vector field $\phi_{1} \nu$ generated by the first eigenfuntion $\phi_{1}$ on $\gamma$ and extend to $X \in \mathcal{X}_{\tan} (D^{2})$. Then denote $\psi_{t}$ to be a flow generated by $X$ and define $\psi_{t}(\gamma) = \gamma_{t}$. For small $|t|$, we can expand the geodesic curvature
\begin{equation}
    \kappa_{\gamma_{t}} = \lambda_{1} \phi_{1} |t|  + O(t^{2})
\end{equation}
toward $\gamma$. Moreover, by the definition of $C_{\gamma}$ in (61) and our setting of $X$,
\begin{equation}
    \lim _{t \rightarrow 0} C_{t} = C_{\gamma}.
\end{equation}
By (62)-(63), there exists $\epsilon>0$ such that $\{ \gamma_{t} \}_{t \in [-\epsilon,\epsilon]}$ is a foliation satisfying the conditions (2)-(3) in the statement. The Robin boundary condition in the equation (4) gives the orthogonality at boundary of mean convex foliation.
\end{proof}
\end{prop}
Now we can apply the ideas to construct the squeezing homotopy via flow by Theorem 4.6 and its pullback homotopy in the general geodesic case. We obtain the following lemma:
\begin{lem}
Let $(D^{2},\partial D^{2},g)$ be a Riemannian $2$-disk with a strictly convex boundary endowed with a bumpy metric. Let $\gamma$ be a free boundary embedded geodesic, and $X$ be a simplicial complex with finite dimension $k$. There exists $\delta_{0} = \delta_{0}((D^{2},\partial D^{2},g), \gamma)>0$ with the following property:

For $0<\delta<\delta_{0}$, if $\Phi : X \rightarrow IV_{1}(D^{2})$ is a continuous map in the smooth topology such that
\begin{equation*}
    \sup \{ F(\Phi(x), \gamma) : x \in X \} < \delta,
\end{equation*}
then there is a homotopy $H:[0,1] \times X \rightarrow IV_{1}(D^{2})$ such that $H(0,x) = \Phi(x)$ and $H(1,x) = \gamma$ so that $\Phi$ is nullhomotopic.
\end{lem}
Now we discuss the modification in the quantitative $F$-distance estimate along the squeezing map. We only need to modify the arguments of Lemma 6.2 in \cite{K1} which proves the upper bound of $F$-distance along the squeezing homotopy in Lemma 5.4 of \cite{K1} when $\gamma$ is a strictly stable geodesic with ambient negative curvature. We rewrite the free boundary version of Lemma 6.1 in \cite{K1} first.
\begin{lem} [Lemma 6.1, \cite{K1}]
 Suppose $\gamma$ is a strictly stable geodesic on Riemannian $2$-disk $(D^{2},\partial D^{2},g)$ with a strictly convex boundary which achieves negative ambient Gaussian curvature in $N_{h}(\gamma)$. Then there exists $C=C(|\gamma|)>0$ satisfying the following property: For $0<\epsilon<h^{2}$, if an embedded curve $\alpha$ whose each boundary point is on the each component of $\partial D^{2} \cap N_{h}(\gamma)$ satisfies $|\alpha|<|\gamma|+\epsilon$, then
\begin{equation*}
F(\alpha,\gamma)< C(|\gamma|)(h+ \sqrt{\epsilon}).    
\end{equation*}
\end{lem}
For a given free boundary mean convex foliation $\{ \gamma_{t} \}_{t \in [-z,z]}$, we denote $\gamma_{\le z'}$ by
\begin{equation*}
    \gamma_{\le z'} := \{ x \in D^{2} | x \in \gamma_{t} \text{ for some } |t| \le z' \}.
\end{equation*}
\begin{lem} Let $\gamma$ be a strictly stable geodesic on $(D^{2},\partial D^{2},g)$ with a strictly convex boundary, Gaussian curvature $K(z)<0$ for $z \in N_{h}(\gamma)$, and $X$ be a $k$-dimensional simplicial complex. Then there exists $\epsilon_{0}= \epsilon_{0}(\gamma)>0$ and $C = C(\gamma)>0$ satisfying the following property: For $0<\epsilon<\epsilon_{0}$, if $\Phi:X \rightarrow IV_{1}(D^{2})$ is a continuous map in the smooth topology such that
\begin{equation}
    \sup \{ F(\Phi(x), \gamma) : x \in X \} < \epsilon,
\end{equation}
then there is a homotopy $H:[0,1] \times X \rightarrow IV_{1}(D^{2})$ such that $H(0,x) = \Phi(x)$, $H(1,x) = \gamma$ and the following $F$-distance estimate holds along $H$:
\begin{equation}
    \sup \{ F(H(t,x),\gamma) : x \in X \text{ and } t \in [0,1] \} < C(\gamma) \sqrt{\epsilon}.
\end{equation}
\begin{proof} We consider the $N_{h_{0}}(\gamma)$ can be foliated by a free boundary mean convex foliation $\{ \gamma_{t} \}_{t \in [-z_{0},z_{0}]}$ of $\gamma$ thanks to Proposition 6.1 and suppose the ambient Gaussian curvature in $N_{h_{0}}(\gamma)$ is negative. Take $\epsilon_{0} = h_{0}^{2}/10$. Then by Lemma 3.1 of \cite{K1}, Hausdorff distance between $\Phi(x)$ and $\gamma$ is smaller than $h_{0}$ and $\Phi(x)$ is supported in $N_{h_{0}}(\gamma)$.

Now we assume the condition (64). By Lemma 3.1 in \cite{K1}, $\Phi(x)$ is supported in $N_{\sqrt{10 \epsilon}}(\gamma)$. We take $z \in (0,z_{0})$ by $z:= \inf \{z' \in (0,z) \,|\, N_{\sqrt{10 \epsilon}}(\gamma) \subseteq \gamma_{\le z'}\}$. Note that $\gamma_{\le z} \subseteq N_{2C_{\gamma} \sqrt{10 \epsilon}}(\gamma)$ by Proposition 6.2. Moreover, $|\Phi(x)|< |\gamma|+ \epsilon$ by the definition of the $F$-distance.

We consider the squeezing homotopy $H$ which is the composition of free boundary curve shortening flow and squeezing map of graphical curves in the smaller scale. The length is monotonically decreasing along the free boundary curve shortening flow part in Lemma 6.3 and we still can apply the length bound of (40) and (41) in \cite{K1} for the length bound over the squeezing map. Saying again,
\begin{equation}
    |H(t,x)| <|\gamma| + \epsilon
\end{equation} for $x \in X$ and $t \in [0,1]$. 

Recall $supp(\Phi(x)) \subset N_{\sqrt{10 \epsilon}} \subseteq \gamma_{\le z}$. We apply the avoidance principle of the free boundary curve shortening flow between $\gamma_{z}$ and $\Phi(x)$. Then we have
\begin{equation}
  supp(H(t,x))\subset \gamma_{\le z} \subseteq N_{2C_{\gamma} \sqrt{10 \epsilon}}(\gamma).
\end{equation}
By applying Lemma 6.4 together with (66) and (67), we obtain the $F$-distance upper bound (65) along the homotopy. 
\end{proof}
\end{lem}
We denote $W_{L_{i},j}$ to be the set of the elements in $W_{L_{i}}$ whose support has Morse index less than or equal to $j$. We also have the pulling-tight procedure toward free boundary embedded geodesics with Morse Index upper bound, which is the free boundary version of Theorem 2.3 of \cite{K1}:
\begin{thm}
Let $(D^{2},\partial D^{2},g)$ be a Riemannian $2$-disk with a strictly convex boundary endowed with a bumpy metric. For any minimizing sequence $\{ \Phi_{j} \}$ which is an $i$-sweepout, there is a deformed minimizing sequence $\{\hat{ \Phi}_{j} \}$ of $\{ \Phi_{j} \}$ satisfying the following property. For any small $s>0$, there is some $0<a<L_{i}$ satisfying
    \begin{equation}
        \{ \hat{\Phi}_{j}(x) \in IV_{1}(D^{2}) : |\hat{\Phi}_{j}(x)| \ge L_{i}-a \} \in \bigcup_{\gamma \in  \Lambda(\{ \Phi_{j} \}) \cap W_{L_{i},i}} B^{F}_{s}(\gamma)
        \end{equation}
for all sufficiently large $j$. Moreover, the multiplicity of geodesics in the critical set is $1$.      \end{thm} 

By applying the arguments in the remaining parts of Section 6 and Section 7 in \cite{K1}, we obtain the following Morse Index bound of free boundary embedded geodesic obtained by smooth min-max construction and free boundary curve shortening flow on surfaces.
\begin{thm}
Suppose $(D^{2},\partial D^{2}, g)$ is a Riemannian $2$-disc with convex boundary endowed with a bumpy metric. Then for each $k=1,2$, there exists a free boundary embedded geodesic $\gamma_{k}$ with
\begin{equation*}
    index(\gamma_{k}) = k
\end{equation*}
and these two geodesics satisfy $|\gamma_{1}| <|\gamma_{2}|$.
\end{thm}
\begin{cor}
For a $2$-Riemannian disk $(D^{2},\partial D^{2},g)$ with a convex boundary, for $k=1,2$, there exists a free boundary embedded geodesic $\gamma_{k}$ with
\begin{equation*}
    index(\gamma_{k}) \le k \le index(\gamma_{k})+ nullity (\gamma_{k}).
\end{equation*}
\end{cor}

\end{document}